\documentclass[12pt]{amsart}

\usepackage{amsmath,amsthm,amsfonts,amssymb}
\usepackage{amssymb, amsthm, amsmath}
\usepackage[english]{babel}
\usepackage{amsfonts}

\newcommand{\ra}{\rightarrow}

\newcommand{\ot}{\otimes}
\newcommand{\mtc}{\mathcal}
\newcommand{\lam}{\lambda}
\newcommand{\Lam}{\Lambda}

\newcommand{\al}{\alpha}
\newcommand{\eps}{\epsilon}
\newcommand{\bn}{\begin}
\newcommand{\en}{\end}
\newcommand{\sub}{\subsection}
\newcommand{\mc}{\mathcal}
\newcommand{\D}{\Delta}

\numberwithin{equation}{section}
\newtheorem{lemma}[equation]{Lemma}
\newtheorem{thm}[equation]{Theorem}
\newtheorem{prop}[equation]{Proposition}

\newtheorem{cor}[equation]{Corollary}
\newtheorem{rem}[equation]{Remark}

\newcommand{\lbd}{\Lambda}
\newcommand{\dw}{\downarrow}
\newcommand{\uw}{\uparrow}

\newcommand{\ch}{\chi}
\newcommand{\mtr}{\mathrm}
\numberwithin{equation}{section}
\title[Semisimple Hopf algebras]
{Coset decomposition for semisimple Hopf algebras}

\newcommand{\AMSclasifR}{16W35, \;16W40}

\author{Sebastian  Burciu}
\address{Inst.\ of Math.\ ``Simion Stoilow" of the Romanian Academy
P.O. Box 1-764, RO-014700, Bucharest, Romania\\ smburciu@syr.edu}

\date{December 11, 2007}
\begin{document}
\thanks{MSC (2000): \AMSclasifR}\thanks{The research was supported by CEx05-D11-11/04.10.05.}
\begin{abstract} The notion of double coset for semisimple finite dimensional Hopf algebras is introduced. This is done by considering an equivalence relation on the set of irreducible characters of the dual Hopf algebra. As an application formulae for the restriction of the irreducible characters to normal Hopf subalgebras are given.
\end{abstract}
\maketitle
\section*{Introduction}
In this paper we introduce a notion of double coset for semisimple finite dimensional Hopf algebras, similar to the one for groups. This is achieved considering an equivalence relation on the set of irreducible characters of the dual Hopf algebra.  The equivalence relation that we define generalizes the equivalence  relation introduced in \cite{NR}. Using Frobenius-Perron theory for nonnegative Hopf algebras the results from \cite{NR} are generalized and proved in a simpler manner.

The paper is organized as follows.
In the first Section  we recall some basic results for finite dimensional semisimple Hopf algebras that we need in the other sections.

Section \ref{dcsf} introduce the equivalence relation on the set of irreducible characters of the dual Hopf algebra and it proves the coset decomposition. Using this coset decomposition in the next section we prove a result concerning the restriction of a module to a normal Hopf subalgebra. A formula for the induction from a normal Hopf subalgebra is also described. In the situation of a unique double coset a formula equivalent to the Mackey decomposition formula for groups is described.

Section \ref{dr} considers one of the above equivalence relations but for the dual Hopf algebra. In the situation of normal Hopf subalgebras this relation can be written in terms of the restriction of the characters to normal Hopf subalgebras. Some results similar to those in group theory are proved.

The next sections studies the restriction functor from a semisimple Hopf algebra to a normal Hopf subalgebra. We define a notion of conjugate module similar to the one for modules over normal subgroups of a group. Some results from group theory hold in this more general setting. In particular we show that the induced module restricted back to the original normal Hopf subalgebra has as irreducible constituents the constituents of all the conjugate modules.

Algebras and coalgebras are defined over the algebraically closed ground field $k=\mathbb{C}$. For a vector space $V$ over $k$  by $|V|$ is denoted the dimension $\mtr{dim}_kV$. The comultiplication, counit and antipode of a
Hopf algebra are  denoted by $\Delta$, $\epsilon$ and
$S$, respectively. We use Sweedler's notation $\D(x)=\sum x_1\ot x_2$ for all $x\in H$.  All the other notations are those used in
\cite{Montg}.
\section{Preliminaries}\label{prel}
Throughout this paper $H$ denotes a finite dimensional semisimple Hopf algebra over $K=\mathbb{C}$. It follows that $H$ is also cosemisimple \cite{Lard}. If $K$ is a Hopf subalgebra of $H$ then $K$ is also semisimple and cosemisimple Hopf algebra \cite{Montg}.

Denote by $\mtr{Irr}(H)$ the set of irreducible characters of $H$ and $C(H)$ the character ring  of $H$.  Then $C(H)$ is a semisimple subalgebra of $H^*$ \cite{Z} and $C(H)=\mtr{Cocom}(H^*)$, the space of cocommutative elements of $H^*$. By duality, the character ring of $H^*$ is a semisimple subalgebra of $H$ and under this identification it follows that $C(H^*)=\mtr{Cocom}(H)$.

If $M$ is an $H$-module with character $\chi$ then $M^*$ is also an $H$-module with character $\chi^*=\chi \circ S$. \\This induces an involution $``\;^*\;":C(H)\ra C(H)$ on $C(H)$.

For a finite dimensional semisimple Hopf algebra $H$ use the notation $\Lam_H \in H$ for the integral of $H$ with $\eps(\Lam_H)=|H|$ and $t_H \in H^*$ for the  integral of $H^*$ with $t_H(1)=|H|$. It follows from \cite{Montg} that the regular character of $H$ is given by the formula \begin{equation*}\label{f1}t_{ _H}=\sum_{\ch \in \mtr{Irr}(H)}\chi(1)\chi\end{equation*} The dual formula is
\begin{equation*}\label{f2}\Lambda_{ _H}=\sum_{d \in \mtr{Irr}(H^*)}\eps(d)d\end{equation*}

If $W \in H^*$-mod then $W$ becomes a right $H$-comodule via $\rho :W \ra W \ot H$ given by $\rho(w)=\sum w_{0}\ot w_1$ if and only if $fw=\sum f(w_1)w_0$ for all $w \in W$ and $f \in H^*$.

Let $W$ be a left $H^*$-module. Then $W$ is a right $H$-comodule and  one can associate to it a subcoalgebra of $H$ denoted by $C_{_W}$ \cite{Lar}. If $W$ is simple and $q=|W|$ then $|C_{_W}|=q^2$ and it is a matrix coalgebra. It has a basis $\{x_{ij}\}_{1\leq i,j \leq q}$  such that $\D(x_{ij})=\sum_{l=0}^qx_{il}\ot x_{lj}$ for all $1\leq i,j \leq q$. Moreover $W \cong k<x_{1i}|\;1\leq i \leq q>$ as right $H$-comodule where $\rho(x_{1i})=\D(x_{1i})=\sum_{l=0}^qx_{1l}\ot x_{li}$ for all $1\leq i \leq q$. The character of $W$ as left $H^*$-module is $d \in C(H^*)\subset H$ and it is given by $d=\sum_{i=1}^qx_{ii}$. Then $\eps(d)=q$ and the simple subcoalgebra $C_{ _W}$ is sometimes denoted by $C_d$.

If $M$ and $N$ are right $H$-comodules the the tensor product $M \ot N$ is also a right $H$-comodule. The associated coalgebra of $M \ot N$ is is $CD$ where $C$ and $D$ are the associated subcoalgebras of $M$ and $N$ respectively (see \cite{NR'}).

\section{Double coset formula for cosemisimple Hopf algebras}\label{dcsf}
In this section let $H$ be a semisimple finite dimensional Hopf algebra as before and $K$ and $L$ be two Hopf subalgebras. Then $H$ can be decomposed as sum of $K-L$ bimodules which are free both as left $K$-modules and right $L$-modules and are analogues of double cosets in group theory. To the end of this section we give an application in the situation of a unique double coset.

There is a bilinear form $m : C(H^*)\ot C(H^*) \ra k$ defined as follows (see \cite{NR}). If $M$ and $N$ are two $H$-comodules with characters $c$ and $d$ then  $m(c,\;d)$ is defined as $\mtr{dim}_k\mtr{Hom}^H(M,\;N)$. The following properties of $m$ (see \cite{NR}) will be used later:
$$m(x, \;yz)=m(y^*,\;zx^*)=m(z^*,\;x^*y)\;\;\text{and}\;\;$$$$m(x,y)=m(y,\;x)=m(y^*,\;x^*)$$ for all $x,y,z \in C(H^*)$.

Let $H$ be a finite dimensional cosemisimple Hopf algebra and $K$,
$L$ be two Hopf subalgebras of $H$. We define an equivalence
relation $r_{ _{K,\;L}}^H$ on the set of simple coalgebras of $H$ as following: $C
\sim D$ if $ C \subset \mathrm{KDL}$.

Since the set of simple subcoalgebras is in bijection with $\mtr{Irr}(H^*)$ the above relation
in terms of $H^*$-characters becomes the following: $c \sim d$ if $m(c\;,\lbd_{ _K}d\lbd_{ _L}) > 0$
where $\lbd_{ _K}$ and $\lbd_{ _L}$ are the integrals of $K$ and $L$ with $\eps(\Lam_{ _K})=|K|$ and $\eps(\Lam_{ _L})=|L|$ and $c,\;d \in \mtr{Irr}(H^*)$.

It is easy to see that $\sim$ is an equivalence
relation. Clearly $c \sim c$ for any $c \in \mtr{Irr}(H^*) $ since both $\Lam_{ _K}$ and $\Lam_{ _L}$ contain the trivial character.

Using the above properties of the bilinear form $m$, one has that if $c \sim d$ then $m(d,\; \Lam_{ _K}c\Lam_{ _L})=m(\Lam_{ _K}^*,\;c\Lam_{ _L}d^*)=m(c^*,\;\Lam_{ _L}d^*\Lam_{ _K})=m(c,\;\Lam_{ _K}^*d\Lam_{ _L}^*)=m(c,\;\Lam_{ _K}d\Lam_{ _L})$ since $\Lam_{ _K}^*=\Lam_{ _K}$ and $\Lam_{ _L}^*=\Lam_{ _L}$. Thus $d \sim c$.

The transitivity can be easier seen that holds in terms of simple subcoalgebras. Suppose that $c\sim d$ and $d \sim e$ and $c,\;d,$ and $e$ are three irreducible characters associated to the simple subcoalgebras $C,\;D$ and $E$ respectively. Then $C\subset KDL$ and $D \subset KEL$. The last relation implies that $KDL \subset K^2EL^2=KEL$. Thus $C \subset KEL$ and $c \sim e$.

If $\mathcal{C}_1,\mathcal{C}_2,\cdots \mathcal{C}_l$ are
the equivalence classes of $r_{ _{K,\;L}}^H$ on $\mtr{Irr}(H^*)$ then let \begin{equation}\label{pregnv}a_i=\sum\limits_{d \in
\mtc{C}_i}\eps(d)d\end{equation} for $0 \leq i \leq l$.

For any character $d \in C(H^*)$ let $L_d$ and $R_d$ be the left and
right multiplication with $d$ on $C(H^*)$.
\begin{prop} \label{eigenvv}With the above notations it follows that $a_i$ are eigenvectors of the operator $T=L_{ _{\Lam_{ _K}}}\circ R_{ _{\Lam_{ _L}}}$ on $C(H^*)$ corresponding to the eigenvalue $|K||L|$.
\end{prop}
\bn{proof}  Definition of $r_{ _{K,\;L}}^H$ implies that $d \sim d'$ if and only if $m(d',\;T(d)) >0$. Since $T(\Lam_{ _H})=|H||K|\Lam_{ _H}$ and $\Lam_{ _H}=\sum_{i=1}^la_i$ it follows that $T(a_i)=|K||L|a_i$ for all $1 \leq i \leq \l$.
\end{proof}
\begin{rem}Let $C_1$ and $C_2$ be two subcoalgebras of $H$ and $K=\sum_{n\geq 0}C_1^n$ and $L=\sum_{n \geq 0}C_2^n$ be the two Hopf subalgebras of $H$ generated by them \cite{NR}. The above equivalence relation can be written in terms of characters as follows: $c \sim d$ if
$m(c,\;c_1^ndc_2^m) > 0$ for some natural numbers $m,n\geq 0$.
\end{rem}

In the sequel, we use the Frobenius-Perron theorem for matrices with
nonnegative entries (see \cite{F}). If $A$ is such a matrix then $A$
has a positive eigenvalue $\lambda$ which has the biggest absolute
value among all the other eigenvalues of $A$. The eigenspace
corresponding to $\lambda$ has a unique vector with all entries
positive. $\lambda$ is called the principal value of $A$ and the
corresponding positive vector is called the principal vector of $A$. Also the eigenspace of $A$ corresponding to $\lam$ is called the principal eigenspace of $A$.

The following result is also needed:

\begin{prop}(\cite{F}, Proposition 5.)\label{transpose}
Let $A$ be a matrix with nonnegative entries such that $A$ and $A^t$
have the same principal eigenvalue and the same principal vector.
Then after a permutation of the rows  and the same permutation of the columns $A$ can be decomposed in diagonal blocks $A=\{A_1,\; A_2,\;\cdots,\;A_s\}$
with each block an indecomposable matrix.
\end{prop}

Recall from \cite{F} that a matrix $A\in\mtr{M}_n(k)$ is called decomposable if the set $I=\{1,2, \cdots , n\}$ can be written as a disjoint union $I=J_1\bigcup J_2$ such that $a_{uv}=0$ whenever $u \in J_1$ and $v \in J_2$. Otherwise the matrix $A$ is called indecomposable.

\begin{thm}\label{mainMack}
Let $H$ be a finite dimensional semisimple Hopf algebra over the
algebraically closed field $k$ and $K$, $L$ be two Hopf subalgebras of
$H$. Consider
the linear operator $T=L_{ _{\Lam_{ _K}}}\circ R_{ _{\Lam_{ _L}}}$ on the character
ring $C(H^*)$ and $[T]$ the matrix associated to $T$ with
respect to the standard basis of $C(H^*)$ given by the
irreducible characters of $H^*$.
\begin{enumerate}
\item
The principal eigenvalue of $[T]$ is $|K||L|$.
\item
The eigenspace corresponding to the eigenvalue  $|K||L|$
has \\ ${(a_i)}_{1 \leq i \leq l}$ as $k$- basis were $a_i$ are defined in \ref{pregnv}.
\end{enumerate}
\end{thm}

\begin{proof}
1) Let $\lambda$ be the biggest eigenvalue of $T$ and $v$ the principal
eigenvector corresponding to $\lambda$. Then $\Lam_{ _K}v\Lam_{ _L}=\lambda v$.
Applying $\eps$ on both sides of this relation it follows that
$|K||L|v=\lam v$. Since $v$ has positive entries it follows that
$\lambda=|K||L|$.

2) It is easy to see that the transpose of the matrix $[T]$ is  also $[T]$. To check that
let $x_1,\;\cdots,\;x_s$ be the basis of $C(H^*)$ given by the irreducible characters of $H^*$
and suppose that $T(x_i)=\sum_{j=1}^st_{ij}x_j$. Thus $t_{ij}=m(x_j,\;\Lam_{ _K}x_i\Lam_{ _L})$
and $t_{ji}=m(x_i,\;\Lam_{ _K}x_j\Lam_{ _L})=m(\Lam_{ _K}^*,\;x_j\Lam_{ _L}x_i^*)= m(x_j^
*,\;\Lam_{ _L}x_i^*\Lam_{ _K})=m(x_j, \Lam_{ _K}^*x_i\Lam_{ _L}^*)=t_{ij}$ since $\Lam_{ _K}^*=S(\Lam_{ _K})=\Lam_{ _K}$ and also $\Lam_{ _L}^*=\Lam_{ _L}$. Proposition \ref{transpose} implies that after a permutation of the rows and the same permutation of the columns  the matrix $[T]$ decomposes 
in diagonal blocks $A=\{A_1,\; A_2,\;\cdots,\;A_s\}$ with each block an indecomposable matrix. This decomposition of $[T]$ in diagonal blocks gives a partition $\mtr{Irr}(H^*)=\cup_{i=1}^s\mtc{A}_i$ on the set of irreducible characters of $H^*$, where each $\mtc{A}_i$ corresponds to the rows (or columns) indexing the block $A_i$.

The eigenspace of $[T]$ corresponding to the eigenvalue $\lambda$ is the sum of the eigenspaces of the
diagonal blocks $A_1, A_2,\cdots A_l$ corresponding to the same value. Since each $A_i$ is an
indecomposable matrix it follows that the eigenspace of $A_i$ corresponding to $\lam$ is one
dimensional (see \cite{F}). If $b_j=\sum_{d \in \mc{A}_j}\eps(d)d$ then as in the proof of Proposition \ref{eigenvv} it can be seen that $b_j$ is eigenvector of $T$ corresponding to the eigenvalue $\lam=|K||L|$. Thus $b_j$ is the unique eigenvector of $A_j$ corresponding to the eigenvalue $|K||L|$. Therefore each $a_i$ is a linear combination of these vectors. But if $d \in \mtc{A}_i$ and $d'\in \mtc{A}_j$ with $i \neq j$  then  $m(d', \; T(d))=0$ and the definition of $r_{ _{K,\;L}}^H$ implies that $d \nsim d'$. This means that $a_i$ is a scalar multiple of some $b_j$ and this defines a bijective correspondence between the diagonal blocks and the equivalence classes of the relation $r_{ _{K,\;L}}^H$.
Thus the eigenspace corresponding to the principal eigenvalue  $|K||L|$
has a $k$- basis given by $a_i$ with $0 \leq i \leq l$.
\end{proof}

\begin{cor}\label{bim}
Let $H$ be a finite dimensional semisimple Hopf algebra and
$K,\;L$ be two Hopf subalgebras of $H$. Then $H$ can be decomposed
as
\begin{equation*}
\label{decomp} H=\bigoplus_{i=1}^lB_i
\end{equation*}
where each $B_i$ is a $(K,L)$- bimodule free as both left $K$-module
and right $L$-module.
\end{cor}
\begin{proof}Consider  the above equivalence relation $r_{ _{K,\;L}}^H$, relative to the Hopf subalgebras $K$ and $L$.
For each equivalence class $\mathcal{C}_i$ let $B_i=\bigoplus_{C \in
\mathcal{C}_i}C$. Then $KB_iL=B_i$ from the definition of the
equivalence relation. It follows that
$B_i=KB_iL \in {}_{K}\mathcal{M}_L^H$ which implies that $B_i$ is free
as left $K$-module and right $L$-module \cite{NZ}.
\end{proof}
The bimodules $B_i$ from the above corollary will be called a double coset for $H$ with respect to $K$ and $L$.

\begin{cor}With the above notations, if $d \in \mathcal{C}_i$ then
\begin{equation}\label{formula}\frac{\Lam_{ _K}}{|K|}d\frac{\Lam_{ _L}}{|L|}=\epsilon(d)\frac{a_i}{\epsilon(a_i)}\end{equation}
\end{cor}
\begin{proof}
One has that $\Lam_{ _K}d\Lam_{ _L}$ is an eigenvector of $T=L_{ _{\Lam_{ _K}}}\circ R_{ _{\Lam_{ _L}}}$ with the maximal eigenvalue $|K||L|$. From Theorem \ref{mainMack} it follows that $\Lam_{ _K}d\Lam_{ _L}$ is a linear combination of the elements $a_j$. But $\Lam_{ _K}d\Lam_{ _L}$ cannot contain any  $a_j$ with $j \neq i$ because all the irreducible characters entering in the decomposition of the product are in $\mathcal{C}_i$. Thus $\Lam_{ _K}d\Lam_{ _L}$ is a scalar multiple of $a_i$ and formula \ref{formula} follows.
\end{proof}

\begin{rem}\label{oneside}
Setting $C_1=k$ in Theorem \ref{mainMack} we obtain Theorem 7
\cite{NR}. The above equivalence relation becomes $c \sim d$ if and
only if $m(c,\;dc_2^m) > 0$ for some natural number $m \geq 0$. The
equivalence class corresponding to the simple coalgebra $k1$
consists of the simple subcoalgebras of the powers $C_2^m$ for $m
\geq 0$. Without loss of generality we may assume that this
equivalence class is $\mathcal{C}_1$. It follows that
\begin{equation*}
\frac{d}{\eps(d)}\frac{a_1}{\eps(a_1)}=\frac{a_i}{\eps(a_i)}
\end{equation*}for any irreducible character $d\in \mathcal{C}_i$.
\end{rem}

Let $K$ be a Hopf subalgebra of $H$ and $s=|H|/|K|$. Then $H$ is free as left $K$-module \cite{NZ}. If $\{a_i\}_{i=1,\;s}$ is a basis of $H$ as left $K$-module then $H=Ka_1\oplus Ka_2\cdots\oplus Ka_s$ as  left $K$ modules. Consider the operator $L_{ _{\Lam_ {_K}}}$ given by left multiplication on $H$. The eigenspace corresponding to the eigenvalue $|K|$ has a basis given by $\Lam_{ _K}a_i$, thus it has dimension $s$. If we restrict the operator $L_{ _{\Lam_ {_K}}}$ on $C(H^*)\subset H$ then Theorem \ref{mainMack} implies that the number of equivalence classes of $r_{ _{K,\;k}}^H$ is equal to the dimension of the eigenspace of $L_{ _{\Lam_ {_K}}}$ corresponding to the eigenvalue $|K|$. Thus the number of the equivalence classes of $r_{ _{K,\;k}}^H$ is always less or equal then the index of $K$ in $H$.

\bn{example}Let $H=kQ_8\#^{\alpha}kC_2$ the 16-dimensional Hopf algebra described in \cite{Kash}. Then $G(H^*)=<x> \times <y>\cong \mathbb{Z}_2\times \mathbb{Z}_2$ and $\mtr{Irr}(H^*)$ is given by the four one dimensional characters $1,x,y,xy$ and $3$ two dimensional characters denoted by $d_1,d_2,d_3$. The algebra structure of $C(H^*)$ is given by:
$$x.d_1=d_3=d_1.x,\;\;x.d_2=d_2=d_2.x,\;\;x.d_3=d_1=d_3.x$$
$$y.d_i=d_i=d_i.y \;\;\text{for all}\;i=1,3$$
$$d_1^2=d_3^2=x+xy+d_2,\;\;d_2^2=1+x+y+xy,\;d_1d_2=d_2d_1=d_1+d_3$$ $$d_1d_3=d_3d_1=1+y+d_2$$
Consider $K=k<x>$ as Hopf subalgebra of $H$ and the equivalence relation $r_{ _{k,\;K}}^H$ on $\mtr{Irr}(H^*)$. The equivalence classes are given by $\{1,x\}$, $\{y,xy\}$, $\{d_2\}$ and $\{d_1,d_3\}$ and the number of them is strictly less than the index of $K$ and $H$. If $C_2 \subset H$ is the coalgebra associated to $d_2$ then the third equivalence class determine in the decomposition \ref{decomp} the free $K$-module $C_2K=C_2$ whose rank is less then the dimension of $C_2$.
\end{example}
\section{More on coset decomposition}\label{mcd}
Let $H$ be a semisimple Hopf algebra and $A$ be a Hopf subalgebra. Define $H//A=H/HA^+$ and let $\pi : H \ra H//A$ be the natural module projection. Since $HA^+$ is a coideal of $H$ it follows that $H//A$ is a coalgebra and $\pi$ is also a coalgebra map.

Let $k$ be the trivial $A$-module via the counit $\eps$. It can be checked that $H//A\cong H\ot_{ _A}k$ as $H$-modules via the map $\hat{h}\ra h\ot_{ _A}1$. Thus $\mtr{dim}_kH//A=\mtr{rank}_{ _A}H$.

If $L$ and $K$ are Hopf subalgebras of $H$ define $LK//K:=LK/LK^+$. $LK$ is a right free $K$-module since $LK \in \mtc{M}_K^H$.
A similar argument to the one above shows that $LK//K\cong LK\ot_{ _K}k$ as left $L$-modules where $k$ is the trivial $K$-module. Thus $\mtr{dim}_kLK//K=\mtr{rank}_{ _K}LK $. It can be checked that $LK^+$ is a coideal in $LK$ thus $LK//K$ has a natural coalgebra structure.

\begin{thm} Let $H$ be a semisimple Hopf algebra and $K,\;L$ be two Hopf subalgebras of $H$. Then $L//L\cap K\cong LK//K$ as coalgebras and left $L$-modules.
\en{thm}

\bn{proof} Define the map $\phi: L \ra LK//K$ by $\phi(l)=\hat{l}$. Then $\phi$ is the composition of $L \hookrightarrow LK \ra LK//K$ and is a coalgebra map as well as a morphism of left $L$-modules. Moreover $\phi$ is surjective since $\widehat{lk}=\eps(k)\hat{l}$ for all $l \in L$ and $k \in K$. Clearly $L(L\cap K)^+ \subset \mtr{ker}(\phi)$ and thus $\phi$ induces a surjective map $\phi: L//L\cap K\ra LK//K$.

Next it will be shown that
\bn{equation}\label{dims}\frac{|L|}{|L\cap K|}=\frac{|LK|}{|K|}\en{equation} which implies that $\phi$ is bijective since both spaces have the same dimension.
Consider on $\mtr{Irr}(H^*)$ the  equivalence relation introduced above and corresponding to the linear operator $L_{ _{\Lam_{ _L}}}\circ R_{ _{\Lam_{ _K}}}$.
Assume without loss of generality that $\mathcal{C}_1$ is the equivalence class of the character $1$ and
put $d=1$ the trivial character, in the formula \ref{formula}. Thus $\frac{\Lam_{ _L}}{|L|}\frac{\Lam_{ _K}}{|K|}=\frac{a_1}{\epsilon(a_1)}$.  But from the definition of $\sim$ it follows that $a_1$ is formed by the characters of the coalgebra $LK$.
On the other hand $\Lam_{ _L}=\sum_{d \in Irr(L^*)}\eps(d)d$ and $\Lam_{ _K}=\sum_{d \in Irr(K^*)}\eps(d)d$
(see \cite{Lar}). Equality \ref{dims} follows counting the multiplicity of the
irreducible character $1$ in $\Lam_{ _K}\Lam_{ _L}$. Using \cite{NR} we know
that $m(1\;, dd')>0$ if and only if $d'=d^*$ in which case $m(1\;,
dd')=1$. Then $m(1, \;\frac{\Lam_{ _L}}{|L|}\frac{\Lam_{ _K}}{|K|})=\frac{1}{|L||K|}\sum_{d \in \mtr{Irr}(L\cap K)}\eps(d)^2=\frac{|L\cap K|}{|L||K|}$ and $m(1, \;\frac{a_1}{\epsilon(a_1)})=\frac{1}{\eps(a_1)}=\frac{1}{|LK|}$.
\end{proof}

\begin{cor} If $K$ and $L$ are Hopf subalgebras of $H$ then $\mtr{rank}_{K}LK=\mtr{rank}_{L\cap K}L$.
\end{cor}
\begin{prop}Let $H$ be a finite dimensional cosemisimple Hopf algebra and $K,\;L$ be two Hopf subalgebras of $H$ such that $KL=LK$.
If $M$ is a $K$-module then
\begin{equation*}
M\uw_K^{LK}\dw_L \cong (M \dw_{L\cap K})\uw^L
\end{equation*}
as left $L$-modules.
\end{prop}
\begin{proof}

For any
$K$-module $M$ one has
$$M\uw^{LK}\dw_L=LK\ot_KM$$
while $$(M\dw_{L\cap K})\uw^L=L\ot_{L\cap K}M$$
The previous Corollary implies that $\mtr{rank}_{_ K}LK=\mtr{rank}_{_ {L\cap K}}L$ thus both modules above have the same dimension.

Define the map $\phi: L\ot_{L\cap K}M \ra LK\ot_KM$ by $\phi(l\ot_{L\cap K} m)=l\ot_K m$ which is the composition of $L\ot_{L\cap K}M \hookrightarrow LK\ot_{L\cap K}M \rightarrow LK\ot_KM$. Clearly $\phi$ is a surjective homomorphism of $L$-modules. Equality of dimensions implies that $\phi$ is an isomorphism.

\end{proof}
If $LK=H$ then
the previous theorem is the generalization of Mackey's theorem decomposition for groups in
the situation of a unique double coset. corollary

\section{A dual relation}\label{dr}

Let $K$ be a normal Hopf subalgebra of $H$ and $L=H//K$. Then the natural projection $\pi:H\ra L$ is a surjective Hopf map and then $\pi^* :L^* \ra H^*$ is an injective Hopf map. We identify $L^*$ with its image $\pi^*(L^*)$ in $H^*$. This is a normal Hopf subalgebra of $H^*$. In this section we will study the equivalence relation $r_{ _{L^*,\;k}}^{H^*}$ on $\mtr{Irr}(H^{**})=\mtr{Irr}(H)$.

The following result was proved in \cite{Bker}.
\begin{prop} \label{resn} Let $K$ be a normal Hopf subalgebra of a finite dimensional semisimple Hopf algebra $H$ and $L=H//K$. If $t_L \in L^*$ is the integral on $L$ with $t_L(1)=|L|$ then $\eps_{ _K}\uw_K^H=t_L$ and $t_L\dw_K^H=\frac{|H|}{|K|}\eps_{ _K}$
\end{prop}

\begin{prop}Let $K$ be a normal Hopf subalgebra of a semisimple Hopf algebra $H$ and $L=H//K$. Consider the equivalence relation $r_{ _{L^*,\;k}}^{H^*}$ on $\mtr{Irr}(H)$. Then $\ch \sim \mu$ if and only if their restriction to $K$ have a common constituent.
\end{prop}
\begin{proof}
The equivalence relation $r_{ _{L^*,\;k}}^{H^*}$ on $\mtr{Irr}(H)$ becomes the following: $\chi \sim \mu$ if and only if $m_{ _H}(\ch,\;t_{ _L}\mu)>0$. On the other hand, applying the previous Proposition it follows that: 
\begin{eqnarray*}  
m_{ _H}(\ch,\;t_{ _L}\mu)& \!=\! & m_{ _H}(t_{ _L}^*,\;\mu\ch^*)=m_{ _H}(t_{ _L},\;\mu\ch^*)
\\ & = & m_{ _H}(\eps \uw_{ _{K}}^H,\;\mu\ch^*)=m_{ _K}(\eps_{ _{K}},\;(\mu\ch^*)\dw_{ _K})
\\ & = & m_{ _K}(\eps_{ _{K}},\;\mu\dw_{ _K}\ch^*\dw_{ _K})=m_{ _K}(\ch\dw_{ _K},\:\mu\dw_{ _K})
\end{eqnarray*}
Thus $\ch \sim \mu$ if and only if their restriction to $K$ have a common constituent.
\end{proof}

\begin{thm}\label{restrres}Let $K$ be a normal Hopf subalgebra of a semisimple Hopf algebra $H$ and $L=H//K$. Consider the equivalence relation $r_{ _{L^*,\;k}}^{H^*}$ on $\mtr{Irr}(H)$. Then $\ch \sim \mu$ if and only if $\frac{1}{\ch(1)}\ch\dw_K=\frac{1}{\mu(1)}\mu\dw_K$
\end{thm}
\begin{proof} Let $\mathcal{C}_1,\mathcal{C}_2,\cdots \mathcal{C}_l$ be
the equivalence classes of $\sim$ on $\mtr{Irr}(H)$ and let \begin{equation}\label{pregnvddual}a_i=\sum\limits_{\ch \in
\mtc{C}_i}\ch(1)\ch\end{equation} for $0 \leq i \leq l$. If $\mathcal{C}_1$ is the equivalence class of the trivial character $\eps$
then the definition of $r_{ _{L^*,\;k}}^{H^*}$ implies that $a_1=t_{ _L}$. Formula from Remark \ref{oneside} becomes

\begin{equation*}
\frac{\ch}{\ch(1)}\frac{t_{ _L}}{|L|}=\frac{a_i}{a_i(1)}
\end{equation*}

for any irreducible character $\ch \in \mathcal{C}_i$.

Restriction to $K$ of the above relation combined with Proposition \ref{resn} gives: \bn{equation}\label{scform}\frac{\ch\dw_{ _K}}{\ch(1)}=\frac{a_i\dw_{ _K}}{a_i(1)}\end{equation}
Thus  $\ch \sim \mu$ if and only if $\frac{\ch\dw_{ _K}}{\ch(1)}=\frac{\mu\dw_{ _K}}{\mu(1)}$.

\sub{Formulae for restriction and induction}\label{indrestr} The previous theorem implies that the restriction of two irreducible $H$-characters
to $K$ either have the same common constituents or they have no common constituents.
Let $t_H$ be the integral on $H$ with $t_H(1)=1$.
One has that $|H |t_{ _H}=\sum_{i=1}^la_i$ as $|H|t_{ _H}$ is the regular character of $H$.
Since $H$ is free as $K$-module \cite{NZ} it follows that the restriction of $|H|t_{ _H}$ to $K$ is the regular character of $K$ multiplied by $|H|/|K|$. Thus $(|H|t_{ _H})\dw_{ _K}=|H|/|K|(|K|t_{ _K})$. But $|K|t_{ _K}=\sum_{\al \in \mtr{Irr}(K)}\al(1)\al$ and Theorem \ref{restrres}
implies that the set of the irreducible characters of $K$ can be partitioned
in disjoint subsets $\mtc{A}_i$ with $1\leq i \leq l$ such that
\bn{equation*} \label{restrform}
a_i\dw_{ _K}=\frac{|H|}{|K|}\sum_{\al\in \mtc{A}_i}\al(1)\al
\end{equation*}
Then if $\ch \in \mtc{C}_i$ formula \ref{scform} implies that
\bn{equation*}
\ch\dw_{ _K}=\frac{\ch(1)}{a_i(1)}\frac{|H|}{|K|}\sum_{\al\in \mtc{A}_i}\al(1)\al
\end{equation*}
Let $|\mtc{A}_i|=\sum_{\al \in \mtc{A}_i}\al^2(1)$. Evaluating at $1$ the above equality one gets $a_i(1)=\frac{|H|}{|K|}|\mtc{A}_i|$. By Frobenius reciprocity the above restriction formula implies that if $\al \in \mtc{A}_i$ then
\bn{equation}\label{indform}\al\uw_{ _K}^H=\frac{\al(1)}{a_i(1)}\frac{|H|}{|K|}\sum_{\ch \in \mtc{C}_i}\ch(1)\ch=\frac{\al(1)}{a_i(1)}\frac{|H|}{|K|}a_i\end{equation}
\end{proof}


\section{Restriction of modules to normal Hopf subalgebras}\label{restr}

Let $G$ a finite group and $H$ a normal subgroup of $G$. If $M$ is an irreducible $H$-module then
\bn{equation*}
M\uw_H^G\dw_H^G={\oplus_{i=1}^s} \;^{g_i}N
\end{equation*}
where $^gN$ is a conjugate module of $M$ and $\{g_i\}_{i=1,s} $ is a set of representatives for the left cosets of $H$ in $G$. For $g \in G$ the $H$-module $\;^gN$  has the same underlying vector space as $N$ and the multiplication with $h \in H$ is given by $h.n=(g^{-1}hg)n$ for all $n \in N$. It is easy to see that $\;^gN \cong \;^{g'}N$ if $gN=g'N$.

Let $K$ be a normal Hopf subalgebra of $H$ and $M$ be an irreducible $K$-module. In this section we will define the notion of conjugate module
to $M$ similar to group situation. If $d\in \mtr{Irr}(H^*)$ we define a conjugate module $\;^dM$. The left cosets of $K$ in $H$ correspond to the equivalence classes of $r_{ _{K,\;k}}^H$. We will show that if $d,\;d'$ are two irreducible characters in the same equivalence class of $r_{ _{K,\;k}}^H$ then the modules $\;^dM$ and $\;^{d'}M$ have the same constituents.  We will show that the irreducible constituents of $M\uw_K^H\dw_K^H$ and $\oplus_{d \in \mtr{Irr}(H^*)}\;^dM$ are the same.

\bn{rem} Since $K$ is a normal Hopf subalgebra it follows that $\Lam_{ _K}$ is a central element of $H$ (see \cite{Mas'}) and by their definition 
$r_{ _{K,\;k}}^H=r_{ _{k,\;K}}^H$. Thus the left cosets are the same with the right cosets in this situation.

\end{rem}
\sub{}\label{conjdef}
Let $K$ be a normal Hopf subalgebra of $H$ and $M$ be a $K$-module. If $W$ is an $H^*$-module then $W\ot M$ becomes a $K$-module with

\bn{equation}\label{def}
 k(w\ot m)=w_0\ot(S(w_1)kw_2)m
\end{equation}

In order to check that $W\ot M$ is a $K$-module one has that
\bn{eqnarray*}
k.(k'.(w\ot m))& \!=\! &k(w_0\ot(S(w_1)k'w_2)m)\\ & = & w_0\ot (S(w_1)kw_2)(S(w_3)k'w_4)\\ & = & w_0\ot(S(w_1)kk'w_2)m\\ & = & (kk').(w\ot m)
\end{eqnarray*} for all $k,k' \in K$, $w \in W$ and $m \in M$.

It can be checked that if $W \cong W'$ as $H^*$-modules then $W\ot M \cong W'\ot M$. Thus for any irreducible character $d \in \mtr{Irr}(H^*)$ associated to a simple $H$-comodule $W$ one can define the $K$-module $\;^dM \cong W\ot M$.
\begin{prop}Let $K$ be a normal Hopf subalgebra of $H$ and $M$ be an irreducible $K$-module with character $\al \in C(K)$. Suppose that $W$ is a simple $H^*$-module with character $d \in \mtr{Irr}(H^*)$. Then the character $\al_{ _d}$ of the $K$-module $\;^d M$ is given by the following formula:
\bn{equation*}
\al_{ _d}(x)=\al(Sd_1xd_2)
\end{equation*}
for all $x \in K$.
\end{prop}
\begin{proof}
Indeed one may suppose that $W\cong k<x_{1i}\;|\;1\leq i \leq q>$ where
$C_{ _d}=k<x_{ij} \;|\; 1\leq i,j\leq q>$ is the coalgebra associated to $W$ and $q=\eps(d)=|W|$. Then formula \ref{def} becomes
$k(x_{1i}\ot w)=\sum_{j,\;l=1}^qx_{1j}\ot (S(x_{jl})kx_{lj})m$. Since $d=\sum_{i=1}^qx_{ii}$ one gets the the formula for the character $\al_{ _d}$.
\end{proof}

For any $d \in \mtr{Irr}(H^*)$ define the linear operator
$c_{ _d}:C(K)\ra C(K)$ which on the basis given by the irreducible characters is given by $c_{ _d}(\al)=\al_{ _d}$ for all $\al \in \mtr{Irr}(K)$.
\bn{rem} From the above formula it can be directly checked that $^{dd'}\al=\;^{d}(\;^{d'}\al)$ for  all $d, d' \in \mtr{Irr}(H^*)$ and $\al \in C(K)$. This shows that $C(K)$ is a left $C(H^*)$-module. Also one can verify that $\;^d(\al^*)=(\;^d\al)^*$.
\end{rem}
\bn{prop}
Let $K$ be a normal Hopf subalgebra of $H$ and $M$ be an irreducible $K$-module with character $\al \in C(K)$. If $d,\;d' \in \mtr{Irr}(H^*)$ lie in the same coset of $r^H_{ _{k,\;K}}$ then $\;^dM$ and $\;^{d'}M$ have the same irreducible constituents. Moreover $\frac{\al_d}{\eps(d)}=\frac{\al_{d'}}{\eps(d')}$
\end{prop}
\bn{proof}
 Consider the equivalence relation $r_{ _{k,\;K}}^H$ from section \ref{dcsf} and $H=\oplus_{i=1}^sB_i$ the decomposition from Corollary \ref{bim}.
Let $\mtc{B}_1,\cdots,\mtc{B}_s$ be the equivalence classes and $b_i$ defined as in \ref{pregnv}. Then formula \ref{oneside} becomes
\bn{equation}
\frac{d}{\eps(d)}\frac{\Lam{ _K}}{|K|}=\frac{b_i}{\eps(b_i)}
\end{equation}
where  $\Lam{ _K}$ is the integral in $K$ with $\eps(\Lam{ _K})=|K|$ and $d \in \mtc{B}_i$. Thus
 \bn{eqnarray*}
\al_{ _{b_i}}(x)\! & = & \! \al(S(b_i)_1x(b_i)_2)=
 \\ & = & \frac{\eps(b_i)}{\eps(d)|K|}\al(S({d\Lam_{ _K})}_1)x(d{\Lam_{ _K}})_2)=
\\ & = & \frac{\eps(b_i)}{\eps(d)|K|}\al(S(({\Lam_{ _K}})_1)S(d_1)xd_2({\Lam_{ _K}})_2)=
\\ & =& \frac{\eps(b_i)}{\eps(d)}\al(S(d_1)xd_2)=
\\ & =& \frac{\eps(b_i)}{\eps(d)}\al_{ _d}(x)
\end{eqnarray*}
for all $d \in \mtc{B}_i$.

This implies that $d\sim d'$ then $\frac{\al_{ _d}}{\eps(d)}=\frac{\al_{ _{d'}}}{\eps(d')}$.
\end{proof}

Let $N$ be an $H$-module and $W$ an $H^*$-module. Then $W\ot N$ becomes an $H$-module such that

\bn{equation}\label{def'}
 h(w\ot m)=w_0\ot(S(w_1)hw_2)m
\end{equation}
It can be checked that $W\ot N\cong N^{|W|}$ as $H$-modules. Indeed the map $$\phi: W\ot N \ra \; _{ _{\eps}}W\ot N, \;\;\;w\ot n \mapsto w_0 \ot w_1n$$  is an isomorphism of $H$-modules where $\;_{ _{\eps}}W$ is considered left $H$-module with the trivial action. Its inverse is given by $w\ot m \mapsto w_0 \ot S(w_1)m$. To check that $\phi$ is an $H$-module map one has that
\bn{eqnarray*}
\phi(h.(w\ot n))\! & = & \! \phi(w_0 \ot (S(w_1)hw_2)m)
\\ & = & w_0 \ot w_1 (S(w_2)hw_3)m
\\ & = & w_0 \ot hw_1m
\\ & = & h.( w_0 \ot w_1m)
\\ & = &h\phi(w\ot m)
\end{eqnarray*}
for all $w \in W$, $m \in M$ and $h \in H$.
\bn{prop}
Let $K$ be a normal Hopf subalgebra of $H$ and $M$ be an irreducible $K$-module with character $\al \in C(K)$. If $d\in \mtr{Irr}(H^*)$ then 
\bn{equation*}
\frac{1}{\eps(d)}\al_{ _d}\uw_{ K}^H=\al\uw_{ _K}^H
\end{equation*}
\end{prop}
\bn{proof} Using the notations from subsection \ref{indrestr} let $\mc{A}_i$ be the subset of $\mtr{Irr}(K)$ which contains $\al$. It is enough to show that the constituents of $\al_{ _d}$ are contained in this set and then the induction formula \ref{indform} from the same subsection can be applied. For this, suppose $N$ is an irreducible $H$-module and
\bn{equation}N\dw_{ _K}=\oplus_{i=1}^sN_i
\end{equation}
where $N_i$ are irreducible $K$-modules. The above result implies that $W\ot N \cong N^{|W|}$ as $H$-modules. Therefore $(W\ot N)\dw_{ _K}={N\dw_{ _K}}^{|W|}$ as $K$-modules.
But $(W\ot N)\dw_{ _K}=\oplus_{i=1}^s(W \ot N_i)$ where each $W \ot N_i$ is a $K$-module by \ref{conjdef}. Thus
\bn{equation}
\oplus_{i=1}^sN_i^{|W|}=\oplus_{i=1}^s(W \ot N_i)
\end{equation}
This shows that if $N_i$ is a constituent of $N\dw_{ _K}$ then $W\ot N_i$ has all the irreducible $K$- constituents among those of $N\dw_{ _K}$.
Formula \ref{indform} applied for each irreducible constituent of $\al_{ _d}$ gives  that
\bn{equation}
\frac{1}{\eps(d)}\al_{ _d}\uw_{ K}^H=\al\uw_{ _K}^H
\end{equation}
for all $\al \in \mtr{Irr}(K)$ and $d \in  \mtr{Irr}(H^*)$.
\end{proof}
\bn{prop}
Let $K$ be a normal Hopf subalgebra of $H$ and $M$ be an irreducible $K$-module. Then $M\uw_K^H\dw_K^H$ and $\oplus_{d \in \mtr{Irr}(H^*)}\;^dM$ have the same irreducible constituents.
\end{prop}
\bn{proof}
Consider the equivalence relation $r_{ _{k,\;K}}^H$ from Section \ref{dcsf} and let $\mtc{B}_1,\cdots,\mtc{B}_s$ be its equivalence classes. Pick an irreducible character $d_i \in \mtc{B}_i$ in each equivalence class of $r_{ _{k,\;K}}^H$ and let $C_i$ be its associated simple coalgebra. Then Corollary \ref{bim} implies that $H=\oplus_{i=1}^sC_iK$. It follows that the induced module $M\uw_{ _K}^H$ is given by
$$M\uw_{ _K}^H=H\ot_{ _K}M=\oplus_{i=1}^sC_iK\ot_{ _K}M$$ Each $C_iK\ot_{ _K}M$ is a $K$-module by left multiplication with elements of $K$ since
$$k.(ck'\ot_Km)=c_1(Sc_2kc_3)k'\ot_Km=c_1\ot_K(Sc_2kc_3)k'm$$ for all $k,k' \in K$, $c \in C_i$ and $m \in M$.
Thus $M\uw_{ _K}^H$ restricted to $K$ is the sum of the $K$-modules $C_iK\ot_{ _K}M$. On the other hand the composition of the canonical maps  $C_i\ot K \hookrightarrow C_iK\ot M \rightarrow C_iK\ot_{_ K}M$ is a surjective morphism of $K$-modules which implies that $C_iK\ot_{_ K}M$ is a homeomorphic image of $\eps(d_i)$ copies of $^{d_i}M$. Therefore the irreducible constituents of $M\uw_K^H\dw_K^H$ are among those of $\oplus_{d \in \mtr{Irr}(H^*)}\;^dM$. In the proof of the previous Proposition we showed the other inclusion. Thus $M\uw_K^H\dw_K^H$ and $\oplus_{d \in \mtr{Irr}(H^*)}\;^dM$ have the same irreducible constituents.
\end{proof}

\bibliographystyle{amsplain}
\bibliography{cosetdecomp}
\end{document}